\documentclass[]{amsart}
\usepackage{amssymb,amsmath}
\usepackage[mathscr]{euscript}

\newcounter{sec}

\newcounter{punct}[sec]

\def\punct{\refstepcounter{punct}{\arabic{sec}.\arabic{punct}.  }}

\newtheorem{theorem}{Theorem}[sec]
\newtheorem{proposition}[theorem]{Proposition}

\newtheorem{lemma}[theorem]{Lemma}

\newtheorem{corollary}[theorem]{Corollary}

\def\COUNTERS{\addtocounter{sec}{1}
              \setcounter{punct}{0}
          \setcounter{equation}{0}
          \setcounter{theorem}{0}
         }
          
          \def\sm{\smallskip}
          
\begin{document}

\newcommand{\supp}{\mathop {\mathrm {supp}}\nolimits}
\newcommand{\rk}{\mathop {\mathrm {rk}}\nolimits}
\newcommand{\Aut}{\mathop {\mathrm {Aut}}\nolimits}
\newcommand{\Out}{\mathop {\mathrm {Out}}\nolimits}
\newcommand{\OO}{\mathop {\mathrm {O}}\nolimits}
\renewcommand{\Re}{\mathop {\mathrm {Re}}\nolimits}
\newcommand{\ch}{\cosh}
\newcommand{\sh}{\sinh}

\def\0{\mathbf 0}

\def\ov{\overline}
\def\wh{\widehat}
\def\wt{\widetilde}

\renewcommand{\rk}{\mathop {\mathrm {rk}}\nolimits}
\renewcommand{\Aut}{\mathop {\mathrm {Aut}}\nolimits}
\renewcommand{\Re}{\mathop {\mathrm {Re}}\nolimits}
\renewcommand{\Im}{\mathop {\mathrm {Im}}\nolimits}
\newcommand{\sgn}{\mathop {\mathrm {sgn}}\nolimits}
\newcommand{\Isoc}{\mathop {\mathrm {Isoc}}\nolimits}
\newcommand{\PIsoc}{\mathop {\mathrm {PIsoc}}\nolimits}

\newcommand{\Sch}{\mathop {\mathrm {Sch}}\nolimits}
\newcommand{\sch}{\mathop {\mathrm {sch}}\nolimits}
\newcommand{\Fr}{\mathop {\mathrm {Fr}}\nolimits}
\newcommand{\Op}{\mathop {\mathrm {Op}}\nolimits}
\newcommand{\Mat}{\mathop {\mathrm {Mat}}\nolimits}
\newcommand{\Exp}{\mathop {\mathrm {Exp}}\nolimits}

\def\bfa{\mathbf a}
\def\bfb{\mathbf b}
\def\bfc{\mathbf c}
\def\bfd{\mathbf d}
\def\bfe{\mathbf e}
\def\bff{\mathbf f}
\def\bfg{\mathbf g}
\def\bfh{\mathbf h}
\def\bfi{\mathbf i}
\def\bfj{\mathbf j}
\def\bfk{\mathbf k}
\def\bfl{\mathbf l}
\def\bfm{\mathbf m}
\def\bfn{\mathbf n}
\def\bfo{\mathbf o}
\def\bfp{\mathbf p}
\def\bfq{\mathbf q}
\def\bfr{\mathbf r}
\def\bfs{\mathbf s}
\def\bft{\mathbf t}
\def\bfu{\mathbf u}
\def\bfv{\mathbf v}
\def\bfw{\mathbf w}
\def\bfx{\mathbf x}
\def\bfy{\mathbf y}
\def\bfz{\mathbf z}

\def\bfA{\mathbf A}
\def\bfB{\mathbf B}
\def\bfC{\mathbf C}
\def\bfD{\mathbf D}
\def\bfE{\mathbf E}
\def\bfF{\mathbf F}
\def\bfG{\mathbf G}
\def\bfH{\mathbf H}
\def\bfI{\mathbf I}
\def\bfJ{\mathbf J}
\def\bfK{\mathbf K}
\def\bfL{\mathbf L}
\def\bfM{\mathbf M}
\def\bfN{\mathbf N}
\def\bfO{\mathbf O}
\def\bfP{\mathbf P}
\def\bfQ{\mathbf Q}
\def\bfR{\mathbf R}
\def\bfS{\mathbf S}
\def\bfT{\mathbf T}
\def\bfU{\mathbf U}
\def\bfV{\mathbf V}
\def\bfW{\mathbf W}
\def\bfX{\mathbf X}
\def\bfY{\mathbf Y}
\def\bfZ{\mathbf Z}
\def\bfT{\mathbf T}

\def\frF{\mathfrak F}
\def\frD{\mathfrak D}
\def\frX{\mathfrak X}
\def\frS{\mathfrak S}
\def\frZ{\mathfrak Z}
\def\frL{\mathfrak L}
\def\frG{\mathfrak G}
\def\frg{\mathfrak g}
\def\frh{\mathfrak h}
\def\frf{\mathfrak f}
\def\frl{\mathfrak l}
\def\frp{\mathfrak p}
\def\frq{\mathfrak q}
\def\frr{\mathfrak r}

\def\bfw{\mathbf w}

\def\R {{\mathbb R }}
 \def\C {{\mathbb C }}
  \def\Z{{\mathbb Z}}
\def\K{{\mathbb K}}
\def\N{{\mathbb N}}
\def\Q{{\mathbb Q}}
\def\A{{\mathbb A}}
\def\U{{\mathbb U}}

\def\T{\mathbb T}
\def\P{\mathbb P}

\def\G{\mathbb G}

\def\cD{\EuScript D}
\def\cK{\EuScript K}
\def\cM{\EuScript M}
\def\cN{\EuScript N}
\def\cP{\EuScript P}
\def\cQ{\EuScript Q}
\def\cR{\EuScript R}
\def\cT{\EuScript T}
\def\cW{\EuScript W}
\def\cY{\EuScript Y}
\def\cF{\EuScript F}
\def\cG{\EuScript G}
\def\cZ{\EuScript Z}
\def\cI{\EuScript I}
\def\cB{\EuScript B}
\def\cA{\EuScript A}
\def\cO{\EuScript O}
\def\cE{\EuScript E}
\def\ex{\EuScript E}

\def\bbA{\mathbb A}
\def\bbB{\mathbb B}
\def\bbD{\mathbb D}
\def\bbE{\mathbb E}
\def\bbF{\mathbb F}
\def\bbG{\mathbb G}
\def\bbI{\mathbb I}
\def\bbJ{\mathbb J}
\def\bbL{\mathbb L}
\def\bbM{\mathbb M}
\def\bbN{\mathbb N}
\def\bbO{\mathbb O}
\def\bbP{\mathbb P}
\def\bbQ{\mathbb Q}
\def\bbS{\mathbb S}
\def\bbT{\mathbb T}
\def\bbU{\mathbb U}
\def\bbV{\mathbb V}
\def\bbW{\mathbb W}
\def\bbX{\mathbb X}
\def\bbY{\mathbb Y}

\def\kappa{\varkappa}
\def\epsilon{\varepsilon}
\def\phi{\varphi}
\def\le{\leqslant}
\def\ge{\geqslant}

\def\B{\mathrm B}

\def\la{\langle}
\def\ra{\rangle}

\def\lambdA{{\boldsymbol{\lambda}}}
\def\alphA{{\boldsymbol{\alpha}}}
\def\betA{{\boldsymbol{\beta}}}
\def\mU{{\boldsymbol{\mu}}}

\def\const{\mathrm{const}}
\def\rem{\mathrm{rem}}
\def\even{\mathrm{even}}
\def\SO{\mathrm{SO}}
\def\SL{\mathrm{SL}}
\def\SU{\mathrm{SU}}
\def\GL{\operatorname{GL}}
\def\End{\operatorname{End}}
\def\Mor{\operatorname{Mor}}
\def\Aut{\operatorname{Aut}}
\def\inv{\operatorname{inv}}
\def\red{\operatorname{red}}
\def\Ind{\operatorname{Ind}}
\def\dom{\operatorname{dom}}
\def\im{\operatorname{im}}
\def\md{\operatorname{\,mod\,}}
\def\St{\operatorname{St}}
\def\Ob{\operatorname{Ob}}
\def\PB{{\operatorname{PB}}}
\def\Tra{\operatorname{Tra}}

\def\ZZ{\mathbb{Z}_{p^\mu}}
\def\F{\mathbb{F}}

\def\cH{\EuScript{H}}
\def\cQ{\EuScript{Q}}
\def\cL{\EuScript{LS}_n}
\def\cX{\EuScript{X}}

\def\Di{\Diamond}
\def\di{\diamond}

\def\fin{\mathrm{fin}}
\def\ThetA{\boldsymbol {\Theta}}

\def\0{\boldsymbol{0}}

\def\FF{\,{\vphantom{F}}_3F_2}
\def\HH{\,\vphantom{H}^{\phantom{\star}}_3 H_3^\star}
\def\Ho{\,\vphantom{H}_2 H_2}

\def\disc{\mathrm{disc}}
\def\cont{\mathrm{cont}}

\def\fan{\vphantom{|^|}}

\def\osigma{\ov\sigma}
\def\ot{\ov t}

\def\Afr{\mathrm{Afr}}
\def\fr{\mathfrak{fr}}
\def\Fr{\mathrm{Fr}}

\def\tri{|\!|\!|}

\begin{center}
\bf\Large
On visualisation of completions of  free groups

\bigskip

\sc Yury A. Neretin%
\footnote{The work is supported by the grant of FWF (Austrian Scientific Funds), PAT5335224.}
\end{center}

{\small The purpose of the note is to visualize the pronilpotent completion of a finitely generated free group
	as a certain subgroup in the free Lie group.

}
\section{Introduction}

\COUNTERS

{\bf\punct Free Lie groups and pronilpotent completions of free groups.}
Denote by $F_n$ the free group with $n$ generators. Let 
$F_n^{(k)}$ be its lower central series. According a general
construction of Malcev \cite{Mal1} (see, also its exposition
in \cite{KM}, Sect.~17),
any finitely generated discrete nilpotent  group  without torsion canonically embeds as a cocompact
lattice to a nilpotent Lie group.
In particular, we can apply this construction to the group
 $F_n(k):=F_n/F_n^{(k)}$ and get a nilpotent Lie group, say $\Fr_n(k)$.
This Lie group corresponds to the real Lie algebra $\fr_n{(k)}$ 
with generators $x_1$, \dots, $x_n$ and relations
 $[x_{i_1},[x_{i_2},\dots,[x_{i_k},x_{i_{k+1}}]\dots]]=0$ for all 
tuples $x_{i_1}$, \dots, $x_{i_{k+1}}$.

 According Malcev \cite{Mal2},
taking
the inverse limit
$$
\Fr_n:=\lim\limits_{\infty \leftarrow k} \Fr_n(k),
$$  
we get an 'infinite-dimensional Lie group', which corresponds to the free Lie algebra
with $n$ generators. 

On the other hand, we have natural surjective homomorphisms $F_n/F_n^{(k+1)}\to F_n/F_n^{(k)}$.
 We consider the inverse limit
\begin{equation}
\cF_n:=\lim\limits_{\infty \leftarrow k} F_n(k),
\label{eq:cF}
\end{equation}
of discrete nilpotent groups $F_n(k)=F_n/F_n^{(k)}$
and get a continual totally disconnected group $\cF_n\subset \Fr_n$
(the {\it pronilpotent completion%
\footnote{Consider a countable discrete group $G$ without torsion. Let $G^{(k)}$
	be its lower central series. Assume that $\cap_k G^{(k)}$ is trivial. The {\it pronilpotent completion} of $G$
is the inverse limit of discrete nilpotent groups $G/G^{(k)}$, see \cite{Bou}, Sect.~II.4. Equivalently,
consider all normal subgroups  $H_\alpha\subset G$, such that  quotient groups $G/H_\alpha$ 
are nilpotent; assume that $\cap_\alpha H_\alpha$ is trivial. If  $H_\alpha\supset H_\beta$,
then we have a canonical map $G/H_\beta\to G/H_\alpha$. So, we can consider the inverse
limit $\lim\limits_{\leftarrow} G/H_\alpha $. It is the same pronilpotent completion,
\newline We vary this definition in two ways. Considering all normal subgroups $H_\alpha$ 
such that quotient groups $G/H_\alpha$ are finite nilpotent groups, we get the {\it pro-finite-nilpotent
completion} $\cG^{\text{fi-ni}}$ of $G$. Considering all normal subgroups such that quotient groups are finite $p$-groups,
we get the {\it pro-$p$-completion} $\cG_p$ of $G$ (by definition, these groups are compact). See \cite{RZ}.
In fact, $\cG^{\text{fi-ni}}$ is isomorphic to the product $\prod_p \cG_p$ over all primes,
see \cite{RZ}, Proposition 2.3.8.
}}
 of $F_n$). 

A desire for transparent descriptions of such objects
arises, in particular, in the context of Kohno  algebras, where the question of producing explicit formulas is natural and important. Recall that Kohno  algebras \cite{Koh} are Lie algebras of Malcev's Lie groups 
corresponding the pure braid groups. They naturally appear in numerous problems of representation theory
and mathematical physics. Quite often actions of Kohno Lie algebras lift to the actions
of the corresponding Lie groups, see \cite{Ner}. On the other hand,
Kohno algebras are normal towers of free Lie algebras (see \cite{Xic}), and the corresponding Lie groups are normal towers of 'free Lie groups' $\Fr_k$ (see \cite{Ner}, Theorem 3.3, Proposition 4.5).

\sm

{\bf \punct The statement of the paper.} Denote by $\Afr_n$
the algebra of formal associative (noncommutative) series in ,,variables''
$\omega_1$, \dots, $\omega_n$.
According Magnus (see, e.g., \cite{LS}, Proposition I.10.1), the multiplicative group
generated by $1+\omega_j\in \Afr_n$ is the free group $F_n$.
 Consider the tensor product
$\Afr_n\otimes \Afr_n$. Consider the homomorphism
of algebras (a coproduct%
\footnote{Our $\delta$ is a coproduct in the sense of Hopf algebras, but formally we 
	do not refer to the structure of a Hopf algebra.})
$$
\delta: \Afr_n\to \Afr_n \otimes \Afr_n 
$$
defined by its values on generators
$$
\delta(\omega_j)=\omega_j\otimes 1+ 1\otimes\omega_j+\omega_j\otimes \omega_j. 
$$

\begin{theorem}
{\rm a)} The group $\Fr_n^\circ$ of elements $g\in \Afr_n$ satisfying the condition
$\delta(g)=g\otimes g$ is isomorphic to $\Fr_n$.

\sm 

{\rm b)} The closure $\cF_n$ of $F_n$
in $\Fr_n^\circ$ is the subgroup consisting of formal series 
with integer coefficients.
\end{theorem}

In fact, classical works by P.~Hall \cite{HalP} and M.~Hall \cite{HalM} give a
 realization
of the pro-nilpotent completion of $F_n$ in terms of Malcev coordinates \cite{Mal1} (see
Magnus, Karras, Solitar, Chap.~5
  Reutenaur \cite{Reu}, Sect.~6.4, but the term pro-nilpotent 'completion' in these books  is not used),  we identify their construction  with
the subgroup $\cF_n$ in $\Fr_n^\circ$ (Corollary \ref{cor:}.b below).

We also present a similar description of the pro-$p$-completions of $F_n$.
 
 \sm 
 
 {\bf \punct The further structure of the paper.} The next section contains preliminaries.
 In Sect.~3 we establish the main results of the paper
 (Theorem \ref{th:} and Corollary \ref{cor:}). Sect.~4 contains some additional details. 
 
 \section{The free Lie group $\Fr_n$. Preliminaries}
 
 \COUNTERS
 
 Let $G$ be a group, by $\{g,h\}$ we denote the group commutator 
 of two elements $g$, $h\in G$:
 $$
 \{g,h\}:=g^{-1}h^{-1} g h.
 $$
 
 By $[x,y]$ we denote the Lie commutator, $[x,y]=xy-yx$. 
 
 \sm
 
 {\bf \punct The free associative algebra.} Consider a metrizable commutative topological ring
 $\cE$ containing
 rational numbers $\Q$, such that the unit of $\Q$ is the unit of $\cE$
 (we assume that the addition and the multiplication in $\cE$ are jointly continuous
 operations).
 Consider the subring $\Z\subset \Q\subset \cE$ and its closure $\ov\Z$ in $\cE$.
 
 The cases, which are actually interesting for us, are:
 
 \sm
 
 --- $\cE$ is the field of reals $\R$ (so $\ov \Z=\Z$);
 
 \sm

 ---  $\cE$ is the field $\Q_p$ of $p$-adic numbers (so $\ov\Z$ is the ring $\Z_p$
 of $p$-adic integers).
 
 \sm
 
 Denote by $\Afr_n(\cE)$ the algebra of associative (noncommutative)
 formal series in letters $\omega_1$, \dots, $\omega_\alpha$.
We denote monomials by 
$$
\omega^\alpha:=\omega_{\alpha_1}\dots \omega_{\alpha_k}
,$$ 
so $\alpha$ ranges in the set $\cA$ of ordered collections
$(\alpha_1,\alpha_2,\dots,\alpha_k)$, where $k=0$, $1$, $2$, $\dots$
and $1\le\alpha_j\le n$.  We write elements of
$\Afr_n(\cE)$ as formal series
$$
x=\sum_{\alpha\in \cA} c_\alpha \omega^\alpha,\qquad c_\alpha\in \cE.
$$
Denote a degree of monomial by 
$$|\alpha|=|(\alpha_1,\alpha_2,\dots,\alpha_k)|:=k.$$

We define the lexicographic order $\prec$ on the set $\cA$ (equivalently,
on the set of monomials $\omega^\alpha$) assuming that the  ``variables'' $\omega_j$ are ordered as
$\omega_1\prec\omega_2\prec \dots\prec \omega_n$. 

We also define another linear order $\lessdot$ on $\cA$:

\sm

--- if $|\alpha|<|\beta|$, then $\omega^\alpha\lessdot\omega^\beta$;

\sm

--- if $|\alpha|=|\beta|$, then $\omega^\alpha\lessdot\omega^\beta$
is
equivalent to $\omega^\alpha\prec\omega^\beta$.

\sm

Denote by $\Afr_n^{[k]}$ the subspace consisting of homogeneous elements
of power $k$ and by $\Afr_n^{[>k]}$ the ideal consisting of
series $\sum_{|\alpha|>k} c_\alpha \omega^\alpha$.
We equip $\Afr_n$ with a topology, assuming that
$\sum c_\alpha^j\omega^\alpha\to \sum c_\alpha\omega^\alpha$
iff we have coefficient-wise convergence $c_\alpha^j\to c_\alpha$
for each $\alpha$. In other words, we consider the space
 $\Afr_n$ as a direct product of a countable number of copies 
of $\cE$ numerated by elements $\alpha\in\cA$
and equip this space with the topology of a direct product.

Let $z_1$, $z_2\in\Afr_n(\cE)$. We write 
$$
z_1=z_2\,\,\bigl(\md \Afr_n^{[>k]}\bigr),
$$ 
if $z_1-z_2\in \Afr_n^{[>k]}$.

For $x\in \Afr_n^{[>0]}$ and $t\in \cE$ we have well defined
\begin{align*}
\exp(x)&:=\sum_{m\ge 0} \frac 1{m!}\,\, x^m\,\,&\in&\,\, 1+\Afr_n^{[>0]};
\\
\ln(1+x)&:=\sum_{m>0} \frac {(-1)^{m+1}}m\,\, x^m\,\,&\in&\,\, \Afr_n^{[>0]};
\\
(1+x)^{t}&:= \sum_{m\ge 0} \frac {t (t-1)\dots(t-m+1)}{m!}\,\, x^m \,\,&\in&\,\, 1+\Afr_n^{[>0]},
\end{align*}
actually, in each homogeneous component $\Afr_n^{[k]}$ we have a finite sum.
These series satisfy the usual identities, as $\ln (\exp(x))=x$, etc.

We need the following trivial lemmas.

\begin{lemma}
	\label{l:integer}
{\rm a)}	Let for $x=\sum_{|\alpha|>0}c_\alpha\omega^\alpha$
all coefficients $c_\alpha$ be integer numbers. Let $t$ also be integer. Then all coefficients 
of the expansion of $(1+x)^t$ are integer.

\sm

{\rm b)} If all coefficients $c_\alpha$ and the exponent $t$ are contained in $\ov\Z$, then all coefficients 
of $(1+x)^t$ are contained in $\ov\Z$.
\end{lemma}

{\sc Proof.} Clearly, for $t\in\Z$
the numbers
$$C_t^m=\frac {t (t-1)\dots(t-m+1)}{m!}$$
are integer. The map $t\mapsto C_t^m$ is continuous on $\cE$. Since $\Z$
is dense in $\ov \Z$, we have $C_t^m\in \ov \Z$ for $t\in\ov \Z$.
\hfill $\square$

\begin{lemma}
\label{l:trivial} Let $k$, $l\ge 1$.
Let $g=1+p+\pi$, where $p\in \Afr_n^k$, $\pi \in \Afr_n^{[>k]}$.
Let $h=1+q+\kappa$, where $q\in \Afr_n^l$, $\kappa \in \Afr_n^{[>l]}$.
Then 
$$
\{g,h\}=1+[p,q]\,\bigl(\md \Afr_n^{[>k+l]}\bigr).
$$
\end{lemma}

{\sc Proof.} Denote $P=p+\pi$, $Q=q+\kappa$. Then
\begin{multline*}
(1+P)^{-1}(1+Q)^{-1}(1+P)(1+Q)=
(1+P)^{-1}\bigl(1+ (1+Q)^{-1}(1+P)(1+Q)\bigr)
=\\=
(1+P)^{-1}\bigl(1+ (1-Q+Q^2-\dots)(1+P)(1+Q)\bigr).
\end{multline*}
We have $Q^2 P$, $Q PQ$, $\dots\in \Afr_n^{[\ge 2l+k]}$,  $P^2Q$, $PQP, \dots \in \Afr_n^{[\ge l+2k]}$
(in both cases, these elements are contained in $\Afr_n^{[>k+l]}$). So, modulo $\Afr_n^{[>k+l]}$ we get
\begin{multline*}
(1+P)^{-1}(1+P+PQ-QP)=
1+PQ-QP + (-P+P^2-\dots)(PQ-QP) =\\=
1+PQ-QP\,
 \bigl(\md \Afr_n^{[>k+l]}\bigr)=
 1+pq-qp\,(\md \Afr_n^{[>k+l]}\bigr)\,.
 \qquad\square\!\!\!\!\vphantom{.}
\end{multline*}

\sm

{\bf \punct The standard coproduct.}
We denote by 
$\Afr_n(\cE)\otimes \Afr_n(\cE)$ the algebra of formal series
$$\sum_{\alpha,\beta\in\cA} h_{\alpha,\beta} \,
\omega^\alpha \otimes \omega^\beta,\qquad h_{\alpha,\beta}\in \cE.$$
The product is defined from the condition
$$
\bigl(\omega^\alpha \otimes \omega^\beta\bigr)
\cdot \bigl(\omega^{\alpha'} \otimes \omega^{\beta'}\bigr)=
\omega^\alpha \omega^{\alpha'} \otimes \omega^\beta \omega^{\beta'}.
$$

The standard {\it coproduct} $\Delta$ in $\Afr_n(\cE)$
(see \cite{Reu}, Subsect. 1.3, \cite{CP}, Definition 2.4.2) is a homomorphism of algebras over $\cE$
$$
\Delta: \Afr_n(\cE)\to \Afr_n(\cE)\otimes \Afr_n(\cE)
$$
defined on generators by 
$$
\Delta(\omega_j)=\omega_j\otimes 1+1\otimes \omega_j.
$$

So, 
$$
\Delta (\omega_j^p)= \sum_{i=0}^p \frac{p!}{i!\,(p-i)!}\,
 \omega^i_j\otimes \omega^{p-i}_j.
$$

We say that an {\it unshuffle} of a monomial $\omega^\alpha$
is a splitting of the factors $\omega_{\alpha_1}$, \dots, $\omega_{\alpha_k}$
into two disjoint subsets, say $\alpha_-$, $\alpha_+$. In particular,
we get two new monomials $\omega^{\alpha_-}$ and $\omega^{\alpha_+}$.
Then
\begin{equation}
\Delta (\omega^\alpha)=\sum_{(\alpha_-,\alpha_+)}
 \omega^{\alpha_-} \otimes \omega^{\alpha_+},
 \label{eq:unsh}
\end{equation}
where the summation is taken over all unshuffles of $\omega^\alpha$.

\sm

{\bf \punct The free Lie algebra and free Lie group.}
Consider the free Lie algebra $\fr_n(\cE)$ generated by $\omega_1$, \dots, $\omega_n$.
For any multiple commutator of elements $\omega_j$ we assign
a polynomial expression in $\omega_1$, \dots, $\omega_n$, 
assigning to any subcommutator $[x,y]$ the associative expression
$xy-yx$. For instance,
$$[\omega_1,\omega_2]\mapsto \omega_1\omega_2-\omega_2\omega_1,\quad
[[\omega_1,\omega_2],\omega_5]\mapsto 
 (\omega_1\omega_2-\omega_2\omega_1)\omega_5-\omega_5(\omega_1\omega_2-\omega_2\omega_1),$$
  etc. Then we get an embedding 
 $\fr_n(\cE)\to \Afr_n$, see, e.g. \cite{Reu}, Theorem 1.4. So we  consider
 $\fr_n(\cE)$ as a subspace in  $\Afr_n(\cE)$. Denote
 by $\ov \fr_n(\cE)$ the closure of $\fr_n(\cE)$ with respect to the topology
 of  $\Afr_n(\cE)$, i.e., the Lie algebra consisting of all formal series
 $$
 z=\sum_{k>0} z_k, \qquad\text{where $z_k\in \fr_n\cap \Afr_n^{[k]}$}.
 $$
 
By the Campbell--Hausdorff formula (see, e.g., \cite{Reu}, Chapter 3), for any $u$, $v\in \ov\fr_n(\cE)$
there exists $z(u,v)\in \ov\fr_n(\cE)$ such that
$$
\exp(u)\,\exp(v)=\exp \bigl(z(u,v)\bigr),
$$
so the set 
$$
\Fr_n(\cE):=\exp\bigl(\ov\fr_n(\cE)\bigr)
$$
is a group with respect to the multiplication in $\Afr_n$.
 The following statements (K.~O.~Friedrichs, R.~Ree) are well known, see, e.g., \cite{Reu},
Theorem 1.4, Theorem 3.2:

\begin{theorem}
	\label{th:fried}
{\rm a)} An element $z\in \Afr_n(\cE)$ is contained in $\ov\fr_n(\cE)$
if and only if
$$
\Delta(z)=z\otimes 1+1\otimes z.
$$

{\rm b)} An element $g\in \Afr_n(\cE)$ is contained in $\Fr_n(\cE)$
if and only if
$$
\Delta(g)=g\otimes g.
$$
\end{theorem}

Keeping in mind this theorem and \eqref{eq:unsh}, 
we see that $g=\sum c_\alpha\omega^\alpha\in\Fr_n(\cE)$
if and only if the coefficients $c_\alpha$ satisfy
 the system of equations
$$
\begin{cases}
c_\alpha c_\beta=\sum_\gamma c_\gamma,
\end{cases}
$$
 where $\gamma$ ranges in the set of all {\it shuffles} (see, e.g., \cite{Reu}) of $\alpha$ and $\beta$.
 This means that we consider all monomials $\omega^\gamma$ equipped with splitting of factors $\omega_{\gamma_j}$
 into two disjoint subsets, elements of the first subset form the product $\omega^\alpha$,  
 elements of the second subset form the product $\omega^\beta$; we take sum of such $c_\gamma$
 with multiplicities (generally, a fixed $\omega^\gamma$ has several admissible splittings).
   For instance,
 $$
 c_{1,2}\,c_{3,4}=c_{1,2,3,4}+c_{1,3,2,4}+c_{1,3,4,2}+c_{3,1,2,4}+c_{3,1,4,2}+c_{3,4,1,2}.
 $$
 
 Notice, that by these equations {\it the group
  $\Fr_n(\cE)$ is closed in   $\Afr_n(\cE)$}.

 \sm 
  
{\bf \punct A twisted  coproduct.}
Consider the automorphism $\gamma$ of $\Afr_n(\cE)$
defined on generators by
$$\gamma(\omega_j):= \ln(1+\omega_j),$$
so $\gamma^{-1}(\omega_j)=\exp(\omega_j)-1$.
 The automorphism $\gamma$
preserves the filtration $\Afr_n^{[>k]}$ and is continuous.
For $t\in \cE$ we have
$$
\gamma\bigl(\exp(t \omega_j)\bigr)=(1+\omega_j)^t.
$$
Conjugating the coproduct $\Delta$ by $\gamma$,  we get a homomorphism
({\it twisted coproduct})
$\delta:\Afr_n(\cE)\mapsto \Afr_n(\cE)\otimes \Afr_n(\cE)$
defined on generators
by
$$
\delta (\omega_j)=\omega_j\otimes 1+1\otimes \omega_j+ \omega_j\otimes \omega_j.
$$
Indeed, we have 
the following chain
\begin{multline*}
\omega_j \xrightarrow{\gamma^{-1}}\, -1+\exp (\omega_j)
\xrightarrow{\Delta}\, -1\otimes 1+ \exp (\omega_j)\otimes \exp (\omega_j)
\,\xrightarrow{\gamma\otimes \gamma}
\\
\xrightarrow{\gamma\otimes \gamma}
 -1\otimes 1+
(1+\omega_j)\otimes (1+\omega_j)=\omega_j\otimes 1+1\otimes \omega_j+ \omega_j\otimes \omega_j.
\end{multline*}

{\sc Remark.} We can consider an arbitrary formal series 
$u(x)=x+\sum_{n\ge 2} h_n x^n$ and an automorphism
$\omega_j\mapsto u(\omega_j)$, see \cite{FPT}, Sect. 5, and get the corresponding
coproduct.
Our statements (Theorem \ref{th:}, Corollary \ref{cor:}, Corollary \ref{cor:p},
Proposition \ref{pr:order}, Proposition \ref{pr:})
  survive  if 
$u(x)=\ln\bigl(1+x+\sum_{j\ge 2} m_j x^j\bigr)$, where
$m_j\in\Z$, but  formulas below require  modifications.
In particular, in Theorem \ref{th:} and Proposition 
\ref{pr:order} we must write
$\Xi_L\bigl(\exp(u(\omega_1)),\dots, \exp (u(\omega_n)  )\bigr)$.
 \hfill $\boxtimes$

\sm

It is easy to see that
\begin{equation*}
\delta(\omega_j^m)
=\sum_{\alpha,\beta: \,0\le\alpha\le m,\, 0\le\beta\le m,\, \alpha+\beta\ge m}
\frac{m!}{(m-\alpha)!\,(m-\beta)!\,(\alpha+\beta-m)!}\,
\omega_j^\alpha\otimes \omega_j^\beta.
\end{equation*}

Consider an arbitrary monomial $\omega^\tau$. We color its letters
green, red, and white. For each such {\it tri-coloring} $\zeta$
we denote by $\omega^{\zeta_+}$ the monomial obtained by omitting
 green letters, by $\omega^{\zeta_-}$ the monomial obtained by
omitting  red letters. 
Then
\begin{equation}
\delta(\omega^\tau)=
\sum_\zeta \omega^{\zeta_+}\otimes \omega^{\zeta_-},
\label{eq:delta-omega}
\end{equation} 
where the summation is taken over the set of all tri-colorings.

Notice that 
\begin{equation*}
\delta(\omega^\alpha)=
\Delta(\omega^\alpha)\,\bigl(\md \Afr_n^{[>|\alpha|]}\bigr).
\label{eq:delta-Delta}
\end{equation*}

\sm

{\bf \punct The group $\Fr_n^\circ(\cE)$.}
Denote  
$$\Fr_n^\circ(\cE):=\gamma(\Fr_n(\cE)).$$
Clearly, we get a topological group isomorphic to
$\Fr_n$. This group is the subset in $\Afr_n(\cE)$
satisfying the equation
$$
\delta(g)=g\otimes g.
$$

Keeping in mind \eqref{eq:delta-omega}, we observe
that $g=\sum c_\alpha \omega^\alpha$ is contained 
in this group iff the coefficients $c_\alpha$ satisfy the following system of quadratic equations:
\begin{equation}
\begin{cases}
c_\alpha c_\beta=\sum_{(\tau,\zeta)} c_\tau,
\end{cases}
\label{eq:system}
\end{equation}
 where $(\tau,\zeta)$ ranges in the set of monomials
 $\omega^\tau$ equipped  with tri-colorings $\zeta$ such that $\omega^{\zeta_+}=
 \omega^\alpha$, $\omega^{\zeta_-}=\omega^\beta$.
 
 \sm


 \begin{lemma}
 \label{l:powers}
 Let $g\in \Fr_n(\cE)$, $t\in \cE$. Then $\gamma(g)^t=\gamma(g^t)$.
 \end{lemma}
 
{\sc Proof.} Denote $h:=\gamma(g)$. First, we assume $\cE=\R$. 
 It
is easy to see that for any $g\in 1+\Afr_n^{[>0]}$
roots $g^{1/m}\in 1+\Afr_n^{[>0]}$ are uniquely defined. Therefore, 
$\gamma$ sends $g^{1/m}$ to $h^{1/m}$. Hence the identity
$\gamma(g^t)=h^t$ holds for any $t\in \Q$, by continuity this is valid for
all $t\in\R$. 

Let us examine an arbitrary $\cE$.
Consider a homogeneous basis $\eta_m$ in $\fr_n(\Q)$.
Let $t_m$, $s\in \cE$. 
We must verify the identity 
$$
\exp\bigl(s\sum\nolimits_{m=0}^\infty t_m\eta_m \bigr)\Bigr|_{\omega_j\mapsto\ln(1+\omega_j)}
=\Bigl(\exp\bigl(\sum\nolimits_{m=0}^\infty t_m\eta_m \bigr)\Bigr|_{\omega_j\mapsto\ln(1+\omega_j)}\Bigr)^s.
$$
In each homogeneous component $\Afr_n^{[k]}$ this is an equality of polynomials in
$t_m$, $s$ with rational coefficients. This equality is valid if $\cE=\R$, and therefore it is valid for any  $\cE$.
\hfill $\square$ 

\section{The closure of the free group}

\COUNTERS

{\bf \punct Lyndon words.} Consider an alphabet of letters $a_1$, \dots, $a_n$. We introduce the lexicographic order $\prec$ on the set of all words.
If a word $u$ is a concatenation of nonempty words $u=vw$, we say that $w$ is a {\it suffix}, and $v$ is 
the {\it corresponding  prefix}. 

A {\it Lyndon word} or {\it Lyndon--Shirshov word} (see a consistent exposition
in  \cite{Lot}, Sect. 5.2) is a word that is $\prec$-less than each its suffix.
Denote the set of all Lyndon words  by $\cL$ (the abbreviation of Lyndon and Shirshov). Denote the length
of $L$ by $|L|$.

 Consider a Lyndon word $u$ of length $>1$ and its longest
suffix $w$ that is a Lyndon word%
\footnote{One letter is a Lyndon word, so the set of Lyndon suffixes is nonempty.}. Then the corresponding
prefix $v$ is also a Lyndon word (see \cite{Lot}, Proposition 5.1.3). We write  $u=(v\circ w)$
with brackets $(\quad)$, after this
we repeat the same reasoning for Lyndon words $v$, $w$, etc. So, we 
get a {\it parenthesization} of the word $u$ (see, e.g., \cite{Lot}),
something like
\begin{equation}
	\label{eq:word}
	a_1a_2a_3a_1a_3a_2a_3\mapsto
((a_1\circ(a_2\circ a_3))\circ((a_1\circ a_3)\circ(a_2\circ a_3))).
\end{equation}

Now we intend to produce three different algebraic expressions from Lyndon words.

\sm

$1^*.$ For an arbitrary  word $M$ (not only a Lyndon word), we define the monomial 
$\omega^{M}$ replacing all letters $a_j$ by the corresponding 
$\omega_j$.

\sm

$2^*.$ For each Lyndon word $L$ we define 
$$
\xi_L(\omega_1,\dots, \omega_n)\in \fr_n(\Q)
$$
by the following rule. We replace each symbol $a_j$ by $\omega_j$ and
each subexpression $(p\circ q)$ in the parenthesization
of $L$ by the commutator $[p,q]$.
For instance, the word $a_1a_2a_3a_1a_3a_2a_3$ (see \eqref{eq:word}) produces
$$
\mapsto ((a_1\circ(a_2\circ a_3))\circ((a_1\circ a_3)\circ(a_2\circ a_3)))
\mapsto
[[\omega_1,[\omega_2, \omega_3]],[[\omega_1, \omega_3],[\omega_2, \omega_3]]]. 
$$
 So, we 
get an element of the homogeneous subspace $\Afr_n^{[|L|]}$,
by the construction all coefficients in this expression are integer.
Then 

\begin{theorem}
Expressions $\xi_L(\omega_1,\dots, \omega_n)$ form a basis of $\fr_n(\Q)$.
\end{theorem}

See, Shirshov \cite{Shi}, see also \cite{Lot}, Theorem 5.3.1.
We need the following statement (it is a special case
of Theorem 5.1, \cite{Reu}).
 
 \begin{lemma}
 \label{l:xi}
\begin{equation}
\xi_L(\omega_1,\dots, \omega_n)=\omega^L+
\sum_{M\succ L} a_{M,L} \omega^M.  
\end{equation}
\end{lemma} 

$3^*.$
Let $L$ be a Lyndon word.
Let $g_1$, \dots, $g_n$ be elements of a group  $G$. We define
$$
\Xi_L(g_1,\dots,g_n)
$$  
in the following way. In $L$ we replace all letters $a_j$ by the corresponding $g_j$ and replace each subexpression $(p\circ q)$
in the parenthesization of $L$ by the group commutator
$\{p,q\}$. For instance, 
$$
\Xi_{a_1a_2a_3a_1a_3a_2a_3}(g_1,g_2,g_3)=
\{\{g_1,\{g_2, g_3\}\},\{\{g_1, g_3\},\{g_2, g_3\}\}\}.
$$

We need the following lemma, see \cite{Reu}, Lemma 6.10
(also, it is a straightforward corollary of  Lemma \ref{l:trivial} above):

\begin{lemma}
\label{l:Xi-xi}
\begin{equation}
\Xi_L(1+\omega_1,\dots, 1+\omega_n)=
1+ \xi_L(\omega_1,\dots, \omega_n)\,\,\bigl(\md \Afr_n^{[>|L|]}\bigr).
\end{equation}
\end{lemma}

{\bf \punct The statements of the paper.}
So, we consider the multiplicative group $\Fr_n^\circ(\cE)\simeq \Fr_n(\cE)$
 consisting of elements
$g\in\Afr_n(\cE)$ satisfying the condition
$\delta(g)=g\otimes g$. 

Denote by $\Fr_n^\circ[\cE,\Z]$ (resp. $\Fr_n^\circ[\cE,\ov\Z]$) the subgroup 
in  $\Fr_n^\circ(\cE)$
consisting of all $g=\sum c_\alpha \omega^\alpha$
such that $c_\alpha$ are contained in $\Z$ (resp. in $\ov\Z$).

\begin{theorem}
\label{th:}
{\rm a)}
Each element of the group $\Fr_n^\circ(\cE)$
can be represented in a unique way as a Malcev type product
\begin{equation} g=
\prod_{L\in\cL}^{(\lessdot)} \Xi_L(1+\omega_1,\dots,1+\omega_n)^{t_L},
\label{eq:product}
\end{equation}
where $t_L$ range in $\cE$ and factors are written
according $\lessdot$-order.

\sm 

{\rm b)} An element $g\in \Fr_n^\circ(\cE)$ is contained 
in $\Fr^\circ_n[\cE,\Z]$ if and only if all $t_L\in \Z$.

\sm

{\rm c)} An element $g\in \Fr_n^\circ(\cE)$ is contained 
in $\Fr^\circ_n[\cE,\ov\Z]$ if and only if all $t_L\in \ov\Z$.
\end{theorem}

Recall that the elements $1+\omega_1$, \dots, $1+\omega_n$ generate
the free group $F_n$. 

\begin{corollary}
\label{cor:}
{\rm a)}
 The closure of $F_n$ in $\Fr_n^\circ(\cE)$
coincides with $\Fr_n^\circ[\cE,\ov \Z]$.

\sm 

{\rm b)} The closure of $F_n$ in $\Fr_n^\circ(\cE)$ consists of all 
products \eqref{eq:product} such that $t_L\in \ov\Z$.
\end{corollary}

Theorem \ref{th:} is proved in the next subsection,
Corollary \ref{cor:} in Subsect. \ref{ss:closure}.

\sm 

{\bf \punct The decomposition of elements of $\Fr^\circ_n(\cE)$.}
Decompose $g\in \Fr^\circ_n(\cE)$ as $g=1+\sum_{k>0} g^{[k]}$, where $g^{[k]}\in\Afr_n^{[k]}$. Consider the first nonzero term, say $g^{[m]}$.
We have 
$$
\delta(1+g^{[m]}+\sum_{k>m} g^{[k]})=
(1+g^{[m]}+\sum_{k>m} g^{[k]})\otimes (1+g^{[m]}+\sum_{k>m} g^{[k]}).
$$ 
Therefore,
$$
\delta(g^{[m]})= g^{[m]}\otimes 1+1\otimes g^{[m]}+
\text{terms of degree $>m$.}
$$
For this reason,
$$
\Delta(g^{[m]})= g^{[m]}\otimes 1+1\otimes g^{[m]}, 
$$
and by Theorem \ref{th:fried} we have $g^{[m]}\in\fr_n(\cE)$. Decompose
$g^{[m]}$ by the basis $\xi_L$,
\begin{equation}
g^{[m]}=\sum_{L\in\cL, |L|=m}
t_L \xi_L(\omega_1,\dots,\omega_m).
\label{eq:t-xi}
\end{equation}
By Lemma \ref{l:Xi-xi},
we have
\begin{multline}
Q:=
\prod_{L\in\cL, |L|=m} \Xi_L(1+\omega_1,\dots,1+\omega_n)^{t_L}
=\\=
1+\sum_{L\in\cL, |L|=m}t_L \xi_L(\omega_1,\dots,\omega_n)\,\,
\bigl(  \md \Afr_n^{[>m]}\bigr)=
1+g^{[m]}\,\, \Bigl(  \mathrm{mod}\,\,\Afr_n^{[>m]}\Bigr)
\label{eq:Q}
\end{multline}
(by Lemma \ref{l:powers}, this product is contained in $\Fr_n^\circ(\cE)$).
Therefore, 
$$
 g Q^{-1}=1\,\, \Bigl(  \mathrm{mod}\,\,\Afr_n^{[>m]}\Bigr).
$$

So, we start by $m=1$ and successively
kill  summands  of $g=1+g^{[1]}+g^{[2]}+\dots$ of degrees 1, 2, 3, \dots.

This proves the statement a) of Theorem \ref{th:}.

\sm 

Proofs of b), c) are identical.  To be
definite, we consider the case $\ov\Z$. 
By Lemma \ref{l:integer}, for $t\in \ov\Z$ 
we have $\Xi_L(\dots)^{t}\in \Fr^\circ[\cE,\ov \Z]$. 
%
%
%
We start the proof as for the statement a). 
All coefficients in the expression
$g^{[m]}$ are contained in $\ov\Z$. We decompose $g^{[m]}$
by the basis $\xi_L=\xi_L(\omega_1,\dots,\omega_n)$, see \eqref{eq:t-xi}. Let $\xi_{L}$, where $L\in \cL$, be the $\prec$-first 
summand
that actually appears in the sum \eqref{eq:t-xi}. 
By Lemma \ref{l:xi}, the $\prec$-first monomial $c_M\omega^M$ that actually
appears in the decomposition of $g^{(m)}$  is $t_L \omega^L$.
Since $t_L$ is a coefficient at monomial, we have $t_L\in\ov \Z$. 
So,
\begin{equation}
g^{[m]}=t_L \xi_L(\omega_1,\dots,\omega_n)+
\sum_{K: K\succ L, |K|=m} t_K \xi_K(\omega_1,\dots,\omega_n).
\label{eq:g-m-xi}
\end{equation}
Now we pass to the element
$$
g':=
g\cdot\Xi_L(1+\omega_1,\dots,1+\omega_n)^{-t_L}.
$$
Then new $(g')^{[m]}$ is
$$
\sum_{K: K\succ L, |L|=m} t_K \xi_K(\omega_1,\dots,\omega_n)
$$
instead of \eqref{eq:g-m-xi}. In the decomposition
$g'=\sum c'_\alpha\omega^\alpha$, all coefficients $c'_\alpha$ remain to be in $\ov\Z$. So we come to the same problem with a shorter sum
\eqref{eq:t-xi}. Repeating the same argument, we 
'kill' $g^{[m]}$.
\hfill $\square$

\sm 

{\bf \punct The closure of the free group.%
\label{ss:closure}} 
By Theorem \ref{th:},  the statements a),  b) of Corollary \ref{cor:} are equivalent.
The group $\Fr_n^\circ[\cE,\ov \Z]$ is closed and contains $F_n$.
%

On the other
hand, for each $L\in \cL$
the subgroup consisting of all
$$\Xi_L(1+\omega_1,\dots,1+\omega_n)^k
,\qquad
\text{where $k$ ranges in $\Z$},
$$
is contained in $F_n$. Decompose this expression as
\begin{equation}
\sum_{m\ge 0} \frac{k(k-1)\dots (k-m+1)}{m!}\,\Bigl(\Xi_L(1+\omega_1,\dots,1+\omega_n)-1 \Bigr)^m.
\label{eq:expr}
\end{equation}
The expression in big brackets has integer coefficients and is contained in
$\Afr_n^{[>0]}$.
Let $t\in\ov \Z$ and a sequence $k_i$ converges to $t$. Then we have coefficient-wise
convergence of  $\Xi_L(\dots)^{k_i}$ to $\Xi_L(\dots)^t$.
Therefore, 
the closure of $F_n$ contains
all elements $\Xi_L(1+\omega_1,\dots,1+\omega_n)^t$, where $t\in\ov \Z$.
Hence  the closure  contains $\Fr_n^\circ[\cE,\ov \Z]$.
\hfill $\square$

\sm

{\bf \punct Pronilpotent and pro-$p$-completions of free groups.}

\sm

a) Let $\cE=\R$. Then 
$\Fr^\circ_n[\cE,\ov \Z]=\Fr^\circ_n[\R,\Z]=\Fr^\circ_n[\R,\ov\Z]$
 is the pronilpotent completion
$\cF_n$ of $F_n$, see \eqref{eq:cF}.

\sm

b) Let $p$ be a prime number. Let $\cE$ be the $p$-adic field $\Q_p$, so $\ov \Z=\Z_p$ is the ring of 
$p$-adic integers. As a corollary, we get the following
statements (it is contained in  Furusho, Enriquez \cite{EF}, Lemma 6.3).

\begin{corollary}
		\label{cor:p}
	 $\Fr_n^\circ[\cE,\ov \Z]=\Fr_n^\circ[\Q_p, \Z_p]$
is the pro-$p$-completion 
of $F_n$.
\end{corollary}

Let us verify this claim.
By definition, the group $\Fr_n^\circ[\Q_p, \Z_p]$ is compact.
Fix  positive integers $\nu$ and $m$ such that $p^m>\nu$.
Consider the subset $U(\nu,p^m)$ in this group consisting of series
$g=1+\sum c_\alpha \omega^\alpha$ such that  for $|\alpha|\le \nu$
we have $c_\alpha\in p^m\Z_p$. It is easy to see that $U(\nu,p^m)$
are open normal subgroups in  $\Fr_n^\circ[\Q_p, \Z_p]$,
and these subgroups form a fundamental system of neighborhoods 
of the unit. So the group $\Fr_n^\circ[\Q_p, \Z_p]$ is profinite,
see, e.g., \cite{RZ}, Theorem 2.1.3.c.  

Next, we notice that 
$$
(1+x_1+x_2+\dots )^{p^m}=1+x_1^{p^m}+x_2^{p^m}+\dots(\md p^m),
$$
 Since $p^m>\nu$,
$$\Xi(1+\omega_1, \dots, 1+\omega_n)^{p^m}=1\Bigl(\md p^m\Afr_n[\Q_p,\Z_p]+\Afr_n^{[> \nu]}[\Q_p,\Z_p]\Bigr).
$$
This easyly implies that a product \eqref{eq:product} is contained
in $U(\nu,p^m)$ if and only if  $t_L\in p^m\Z_p$ for $|L|\le \nu$
and remaining $t_M$ are contained in $\Z_p$. So, the index of our subgroup
is $p^{m\sigma(\nu)}$, where $\sigma(\nu)$ is the number of Lyndon
words of length $\le \nu$. So $\Fr_n^\circ[\Q_p, \Z_p]/U(\nu,p^m)$ is a $p$-group,
and $\Fr_n^\circ[\Q_p, \Z_p]$ is a pro-$p$-group. 








 
 Next, let $\frF_n$ be the pro-$p$-completion of $F_n$. By the universality of $\frF_n$,
 the identical  map $F_n\to F_n$ induces a homomorphism
 $\frF_n\to \Fr^\circ_n[\Q_p,\Z_p]$, see \cite{RZ}, Proposition 3.3.6.
 Since $F_n$ is dense in $\Fr^\circ[\Q_p,\Z_p]$ and $\frF_n$
 is compact, this homomorphism is surjective. Therefore, it is an isomorphism, see \cite{RZ}, Proposition 2.5.2.
 
 \sm

 
c) Consider the adelic sum $\bbA$ of all $p$-adic fields $\Q_p$,
i.e., we  consider the product $\Q_2\times \Q_3\times \Q_5\times\dots$
and the subring $\bbA$ in this product consisting of sequences $(z_2,z_3,z_5,z_7,\dots)$ such that $|z_p|\le 1$ for all but a finite number of primes.
We consider the diagonal embedding $\Z\to \bbA$, i.e., 
$m\mapsto (m,m,m,\dots)$. Then $\ov\Z=\prod_p \Z_p$.

The group 
$$\Fr_n^\circ[\bbA,\ov\Z]\simeq \prod_p \Fr_n^\circ[\Q_p,\Z_p]$$
 is the pro-finite-nilpotent completion of 
$F_n$ (cf. \cite{RZ}, Example 2.1.6, Proposition 2.3.8).

\section{Additional details}

\COUNTERS

{\bf \punct Modifications of Malcev products.}
Consider an arbitrary linear order $\vartriangleleft$ on the set
$\cL$ of all Lyndon words%
\footnote{We admit arbitrary orders. For instance, it can happened
that for any $L\vartriangleleft M$ there is $N$ such that
$L\vartriangleleft N \vartriangleleft M$.}.

\begin{proposition}
\label{pr:order}
	{\rm a)}
Each element of the group $\Fr_n^\circ(\cE)$
can be represented in a unique way as a product
\begin{equation} g=
\prod_{L\in\cL}^{{\vartriangleleft}} \Xi_L(1+u_1,\dots,1+u_n)^{t_L},
\label{eq:product1}
\end{equation}
where $t_L$ range in $\cE$, and factors are ordered according
the order $\vartriangleleft$.

\sm

	{\rm b)} Similar statements hold for the groups $\Fr_n^\circ[\cE,\Z]$ and
	$\Fr_n^\circ[\cE,\ov\Z]$, in these cases, $t_L\in \Z$ and $\ov \Z$ respectively. 
\end{proposition}

\sm 

We must explain a sense of such product.
To evaluate \eqref{eq:product1}, it is sufficient to evaluate
the product modulo $\Afr_n^{[>m]}$ for each $m$. But for all Lyndon words 
$L$ of length $>m$, we have 
$$
\Xi_L(1+u_1,\dots,1+u_n)^{t_L}=1\,\,\bigl(\mathrm{mod}\,\,\Afr_n^{[>m]} \bigr).
$$
So, actually we have a product $g^{(m)}$ of finite number of factors.
 Clearly,
$$g^{(m+1)}=g^{(m)} \,\,\bigl(\md\,\,\Afr_n^{[>m]}\bigr).
$$
Taking limit of $g^{(m)}$ we get our product.  Clearly, it is contained 
in $\Fr_n^\circ(\cE)$.

\sm

{\sc Proof} is the same for all 3 cases.

\sm 

 {\it Uniqueness.} Assume that we have two representations 
of $g$ as a product \eqref{eq:product1}, for a collection $t_L$
and a collection $t_L'$. Take the largest $k$ such that for
$|L|< k$ we have $t_L=t_L'$. By Lemma \ref{l:trivial},

\sm

(*)
{\it if $|M|=k$ and $|L|<k$, then 
  $\Xi_M(1+\omega_1, \dots, 1+\omega_n)$  commutes with
 $\Xi_L(1+\omega_1, \dots, 1+\omega_n)$ modulo $\Afr_n^{[>k]}$.}

\sm 

Therefore, for evaluation of $g^{(k)}$ we can move the factors $\Xi_M(\dots)^{t_M}$ to
the end of the product (leaving the factors $\Xi_L(\dots)^{t_L}$ at their places).
 Since two products \eqref{eq:product1} must coincide
modulo $\Afr_n^{[>k]}$, we have 
$$\prod_{|M|=k} \Xi_M(\dots)^{t_{M}}=\prod_{|M|=k} \Xi_M(\dots)^{t'_{M}}.$$
So $t_M=t'_M$, 
and we come to a contradiction.

\sm

{\it Existence.} Denote by $Q$ the set of all products
\eqref{eq:product1}.
  Consider $g\in \Fr_n^\circ (\cE)$. Consider the largest
$k$ such that there exist $h\in Q$ satisfying 
$g=h \bigl(\md  \Afr_n^{[>k]}\bigr)$. Keeping in mind
the statement (*), we observe that for any collection of $t_M$
the product
$h\prod_{M:\,|M|=k+1}\Xi_M(1+\omega_1,\dots,1+\omega_n)^{t_M}$
equals to some element of $Q$ modulo $\Afr_n^{[>k+1]}$.
But 
$g\bigl(\md \Afr_n^{[>k+1]}\bigr)$ coincide with one of such products.
We again come to a contradiction. 
\hfill $\square$

\sm

{\bf \punct Another parametrization of pronilpotent completions.}

\begin{proposition}
\label{pr:}
{\rm a)} For any collection $a_L\in \cE$, where $L$ ranges in $\cL$,
there is a unique element $g=\sum c_\alpha \omega^\alpha\in \Fr_n^\circ(\cE)$,
such that $c_L=a_L$.

\sm

{\rm b)} The  statement {\rm a)}
 holds for the groups $\Fr_n^\circ[\cE,\Z]$ and
$\Fr_n^\circ[\cE,\ov \Z]$. In these cases, $a_L\in \Z$ and $a_L\in \ov\Z$
respectively.

\sm

{\rm c)} Moreover, each $c_\alpha$, which is not $a_L$, is a polynomial
$P_\alpha(\dots)$ in
$a_M$ with $M\lessdot \alpha$. Each exponent $t_L$ in the product
\eqref{eq:product}  has the form
$t_L=a_L+ q_L$, where $q_L(\dots)$ is a polynomial depending on
$a_M$ with $M\lessdot L$. All these polynomials have rational 
coefficients and take integer values at integer points%
\footnote{But their coefficients are not integer.}. 
\end{proposition}

We need the following lemma:

\begin{lemma}
\label{l:last}
For a Lyndon word $L$ denote by $L_+$ the $\gtrdot$-next Lyndon
word.
 Decompose an element $g=\sum c_\alpha \omega^\alpha\in\Fr^\circ(\cE)$
 into a product
\eqref{eq:product}.
Let $M$ be a Lyndon word. Denote by
$g_M=\sum c^M_\alpha \omega^\alpha$ partial product  for
\eqref{eq:product} over $L\lessdot M$.
Then 

\sm

--- each coefficient $c^M_\alpha$ is a polynomial $P^M_\alpha(\cdot)$
 in 
$c_L$, where $L$ ranges in the set of Lyndon words
$\lessdot M$;

\sm

--- for all $M \gtrdot \alpha$ polynomials
$P^M_\alpha$ coincide;

\sm

--- $t_{M}=c_{M}-P^M_{M}$;

\sm

--- all  polynomials $P^M_\alpha(\cdot)$ have rational 
coefficients and take integer values at integer $a_L$.
\end{lemma}

{\sc Proof of Lemma.} 
%
The stabilization of polynomials $P^M_\alpha$ for $M\gtrdot\alpha$ follows from Lemma \ref{l:xi}.
We prove the remaining statements by induction.
Assume that they are valid for some $g_K$, where $K$ is a Lyndon
word. 
Keeping in mind
 Lemmas \ref{l:Xi-xi}
and \ref{l:xi}, we get
 $$
 g_{K_+}=g_K \cdot \Xi_K(1+\omega_1,\dots,1+\omega_n)^{t_K}
 = g_K\Bigl(1+t_K \omega^K+ \sum_{\gamma\gtrdot K}
 \omega^\gamma\Bigr)   .
 $$ Therefore,
 \begin{equation}
 \label{eq:cct}
c_K= c^{K_+}_{K}=c_{K}^{K}+ t_{K}.
 \end{equation}  
So we can write $g_{K_+}$ as
$$
 g_{K_+}=g_K \cdot
  \Xi_K(1+\omega_1,\dots,1+\omega_n)^{c_{K}-c_{K}^{K}}.
$$  
The coefficients of the second factor are polynomial expressions
in $(c_{K}-P_{K}^K(\cdot))$.
Now all statements for $M=K_+$ became clear. 
 
 The base of the induction is
  $M=\omega_1$ (in this case, $g_{\omega_1}=1$ is  the empty product).
  \hfill $\square$
  
  \sm
  
   {\sc Proof of Proposition \ref{pr:}.}
  The statement c) follows from Lemma 
  \ref{l:last}.  Uniqueness in a) also follows from c).
  Existence follows from \eqref{eq:cct}, since for any 
  $c_{K_+}$ we can choose $t_{K_+}$.
  
  Since polynomials are integer at integer points, we have b).
  \hfill $\square$

\sm 

{\sc Remark.}
 There are other versions of ,,free Lie groups'' \cite{Ner}, \cite{Ner2}.
They are certain dense subgroups in $\Fr_n(\R)$ with stronger topologies, and apparently they are more
interesting than $\Fr_n(\R)$. For such ,,free Lie group'' there arise the following questions:

\sm

--- for which sequences $t_L$ a Malcev product \eqref{eq:product} is contained in this group?

\sm 

--- which are conditions for the coefficients $a_L$ from Proposition \ref{pr:}?

\sm 

--- to describe the closure of $F_n$ in this group.




\tt 

University of Graz,
\\
\phantom{.}
\hfill Department of Mathematics and Scientific computing;

High School of Modern Mathematics MIPT,
\\
\phantom{.}
\hfill 1 Klimentovskiy per., Moscow; 

Moscow State University, MechMath. Dept;

 University of Vienna, Faculty of Mathematics.
 
 \sm

e-mail:yurii.neretin(dog)univie.ac.at

URL: https://www.mat.univie.ac.at/$\sim$neretin/ 

\phantom{URL:} https://imsc.uni-graz.at/neretin/index.html

\end{document}